\PassOptionsToPackage{hyphens}{url}
\documentclass[oneside,twocolumn,10pt]{article}

\usepackage[utf8]{inputenc}
\usepackage[slovak,english]{babel}

\usepackage{graphicx}
\usepackage{biblatex}
\usepackage[dvipsnames]{xcolor} 
\usepackage{amsmath, amsthm, amssymb}
\usepackage{mathtools}
\usepackage{hyperref}
\usepackage{siunitx,booktabs,array} 
\usepackage{subfig}
\usepackage[toc,page]{appendix}

\usepackage[bottom]{footmisc}
\usepackage[total={17cm,25cm}, top=20mm, bottom=20mm, right=20mm, left=20mm, includefoot]{geometry}

\usepackage{tikz}

\title{3-D Shadows of 4-D Algebraic Hypersurfaces in a 4-D Perspective}

\author{
\parbox{0.45\textwidth}{\centering
Jakub {\v{R}}ada\\[1mm]
Charles University\\
Faculty of Mathematics and Physics\\
Mathematical Institute\\
Sokolovská 49/83, 186 00 \\
Prague, Czech Republic\\
~\\[1mm]
rada@karlin.mff.cuni.cz
}
\hspace{0.05\textwidth}
\parbox{0.45\textwidth}{\centering
Michal Zamboj\\[1mm]
Charles University\\
Faculty of Education\\
Department of Mathematics and
Mathematical Education\\
M. Rettigové 4, 116 39\\
Prague, Czech Republic\\[1mm]
michal.zamboj@pedf.cuni.cz
}  
}

\usepackage[hyphens]{url}
\usepackage[hyphenbreaks]{breakurl}

\usepackage{pgfplots}
\pgfplotsset{compat=1.15}
\usepackage{mathrsfs}
\usetikzlibrary{arrows}

\usepackage{xcolor}

\begin{document}

\maketitle  

\begin{abstract}
\noindent 
The paper is focused on the four-dimensional visualization of hypersurfaces represented by implicit equations without their parametrization.  We describe a general method to find shadow boundaries in an arbitrary dimension and apply it in a three- and four-dimensional space. Furthermore, we design a system of polynomial equations to construct occluding contours of algebraic surfaces in a 4-D perspective. The method is presented on a composed 3-D scene and three 4-D cases with gradual complexity. In general, our goal is to improve the understanding of spatial properties in a four-dimensional space.

\end{abstract}

\subsection*{Keywords}
four-dimensional visualization, algebraic hypersurface,  tangent hypercone, shadow, occluding contour, implicit equation, elimination methods  

\vspace*{1.0\baselineskip}


\section{Introduction}

Visualizing shapes embedded in more-dimensional spaces faces several challenges. In many cases, mere projections into three- or less-dimensional spaces contain overlapping parts, making them difficult to understand. One technique that enhances intuition about the properties of shapes and their mutual relations is to visualize shadows cast on themselves and on other objects. Apart from the general case, choosing algebraic hypersurfaces defined by polynomials is often convenient. These are good candidates for visualization using computational methods of algebraic geometry and elimination theory. In this sense, instead of sets of many points and operating with meshes, we can work with implicitly represented hypersurfaces, their projections, contours, intersections, etc. Since algebraic methods preserve many properties of the visualized shapes, they are suitable for precise mathematical visualization. The disadvantage of implicit representation is the computational speed when processing polynomials of higher degrees or adding more variables. 

This paper aims to improve understanding spatial properties in a four-dimensional scene containing algebraic hypersurfaces. To do so, we join theoretical geometric construction and algebraic computational methods and provide concrete examples of visualizations of four-dimensional hypersurfaces and their shadows based on implicit representations. 

In particular, we show visualizations of four-dimensional algebraic hypersurfaces (3-surfaces), their contours (2-surfaces), terminators, and 3-D shadows cast on other 3-surfaces with respect to a point light source. The process consists of two main parts -- central projection of the scene into a 3-D modeling space (usually a virtual 3-D environment in some software, AR, VR, or even a real 3-D model) and construction of shadows from an arbitrary point light source. First, to construct a \mbox{4-D}~perspective image, a 3-surface given by a polynomial is intersected by its first polar (3-surface) with respect to the center of projection, and its (2-surface) contour generator is centrally projected into a modeling 3-space. The second phase is finding the terminator with respect to an arbitrary point light source and its projection to the 3-surface on which the shadow is cast. While geometrically, we describe intersections of surfaces,, algebraically, we need to find a polynomial that solves systems of polynomial equations with several variables (7 to 9 in the 4-D case). This leads to the use of standard computational methods such as finding a Gr{\" o}bner basis or Dixon resultant. Finally, to complete the shadow, especially for 3-surfaces of degrees higher than 2, we need to find the regions in their own shade. Such regions are not algebraically omitted in the elimination procedure; hence, we need to carry out further selection. 

\subsection{Related Work}
The algebraic concepts, in particular, the use of Gr{\" o}bner basis and Dixon resultant for finding solutions of polynomial systems, are described in detail in \cite{Kapur1992}. A similar technique to find the implicit representation of an occluding contour in 3-D through the Dixon resultant was applied in \cite{Khan2007, Khan2014}. In our experiments, we used {\sl Wolfram Mathematica 13} implementations of algorithms for finding Gr{\"o}bner basis --- \mbox{WM-GB} \cite{Wolfram1991} and the Dixon resultant --- WM-Dix \cite{Lichtblau2023}, 
and also Dixon resultant --- Fer-Dix-KSY or improved Fer-Dix-EDF \cite{Dixon2008} implemented in software {\sl Fermat 6.5} \cite{Fermat2008} by Lewis, see \cite{Lewis2008, Lewis2018}. 

In addition to our approach, where we start with surfaces given implicitly, a considerable part of previous research on computational aspects of surfaces deals with implicitization from parametric representation (e.g., \cite{Sederberg1984, Buse2003, Li2004, Lewis2018}). The algebraic derivation of perspective images of surfaces and reconstruction with respect to further applications is shown in \cite{Liu2002}. 

Four-dimensional projections through parameterization or point coordinates are described in \cite{Noll1967, Zacharias2000, Miwa2013}. A 4-D perspective projection was also used to visualize implicitly given surfaces that arise in a complex number plane \cite{Rada2023}.  A descriptive geometry approach for constructing 4-D perspective images of a 3-sphere in 3-D from orthogonal projections is discussed in \cite{Rada2021}. Visualizations of the 4-space, including hypersurfaces, are treated comprehensively in \cite{Banchoff1990}. In \cite{Zhou1991, Hoffmann1991}, the authors created projections of several examples of surfaces in 4-space, examined their properties geometrically and algebraically, and showed various applications. Four-dimensional lighting was used to study the shades of some mathematically interesting 3-surfaces in \cite{Hanson1991}. Interactive manipulation with four-dimensional objects based on their projections or shadows in a hyperplane is elaborated in \cite{Banks1992, Zhang2007} and through a tetrahedral mesh construction in \cite{Chu2009, Cavallo2019}.

\subsection{Contribution}
Our approach emphasizes algebraic methods in four-dimensional visualization. Throughout the paper, we work purely with implicit representations of surfaces without the necessity of their parametrization. Visualizing shadows between 3-surfaces, we discuss their mutual relations in the 4-space. In this way, compared to previous attempts, we offer a more comprehensive perception of complex 4-D scenes projected into 3-D space. We also provide a direct method for visualization of 3-surfaces in a 4-D perspective. After all, the designs of polynomial systems for constructing tangent cones and shadow boundaries as intersections are general for any dimension.

\subsection{Paper organization}
The rest of the article is organized as follows: we start with a 3-D example to describe the algorithm to find the terminator of a 2-surface and its shadow cast on another 2-surface in Section~2. Next, we generalize it into 4D and describe the 4-D perspective from a 4-space into the modeling 3-space. Section~3 is focused on concrete examples. The 3-D scene from the previous explanation is supplemented with technical details. Next, we consecutively examine three 4-D scenes with respect to their geometric and computational complexity. In Section~4, we discuss the critical points of our method and propose further research directions. Section~5 summarizes the results of this paper. The polynomials in their full forms, code in {\sl Mathematica}, computation times, and videos are attached in Appendices.  

\section{Method}
\label{sec:method}
    \subsection{Constructing shadow of an algebraic hypersurface}
    \label{subsection:shadow}
        In the first part, we examine the process of computing a shadow in an arbitrary dimension but visualized in a 3-D case\footnote{The upcoming 3-D visualizations are created in {\sl Wolfram Mathematica 13} from outputs based on implicit equations (or inequalities) formulated in variables $x,y,z$. Therefore, they are also usable as 3-D graphics in an interactive environment.}.

    \begin{figure}[!hb]
        \centering
        \includegraphics[width=\linewidth]{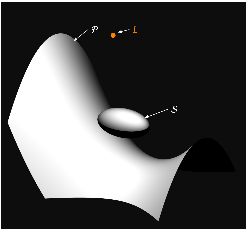}
        \caption{Initial setting of hypersurfaces $\mathcal{S}, \mathcal{P}$, and a point light source $L$.}
        \label{fig:3dset}
        \end{figure}

        \subsubsection{Preliminaries} 
            Let us have a hypersurface $\mathcal{S}$ (Figure~\ref{fig:3dset}), i.e., an \mbox{$(n-1)$-surface} embedded in a real $n$-space $(n\geq 1)$ given by a polynomial equation in $n$ variables  
        \begin{equation}
        \mathcal{S}: \sigma=0.
        \end{equation}

        To find polar hypersurfaces, it will be convenient to work with homogeneous coordinates in the projectively extended real space. 
        
        Let us have $(x'_1,x'_2,\dots,x'_n,x'_0), (y'_1,y'_2,\dots,y'_n,y'_0)\in \mathbb{R}^{n+1}\setminus\{(0,0,\dots,0)\}$. We define the equivalence $(x'_1,x'_2,\dots,x'_n,x'_0)\sim (y'_1,y'_2,\dots,y'_n,y'_0)$, if there exists $\lambda\in\mathbb{R}\setminus\{0\}$ such that $(x'_1,x'_2,\dots,x'_n,x'_0)=(\lambda y'_1,\lambda y'_2,\dots,\lambda y'_n,\lambda y'_0)$. The projective $n$-space is defined as equivalance classes of $\mathbb{R}^{n+1}\setminus\{(0,0,\dots,0)\}$. A point $\overline{P}(p'_1,p'_2,\dots,p'_n,p'_0)$ of the projective space has homogeneous projective coordinates $(p'_1,p'_2,\dots,p'_n,p'_0)$. Additionally, for $p'_0\neq 0$, we can obtain the Cartesian coordinates of $P(p_1,p_2,\dots,p_n)$ by substituting $p_1=\frac{p'_1}{p'_0},p_2=\frac{p'_2}{p'_0}, \dots, p_n=\frac{p'_n}{p'_0}]$, or vice versa. In case $ p'_0=0$, the point $\overline{P}$ represents a point at infinity. Throughout the text, $\overline\sigma, \overline{P},\dots$ denote representations of $\sigma, P, \dots$ in homogeneous coordinates.

        Let $\overline{P}$ be a regular point of $\mathcal{S}$ and assume a polynomial $\overline{\sigma_P}=\overline{P}^T\nabla \overline{\sigma}$. A polar hypersurface $\mathcal{S}_P: \overline{\sigma_P}=0$ is called the first polar of the point $\overline{P}$ with respect to the hypersurface $\mathcal{S}$.
        
        \subsubsection{Terminator} 
        
                 \begin{figure}[!htb]
        \centering
        \includegraphics[width=\linewidth]{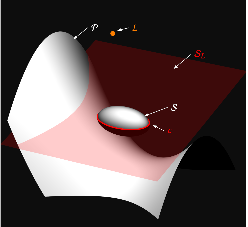}
        \caption{The polar hypersurface $\mathcal{S}_L$ of a hypersurface $\mathcal{S}$ with respect to a point light source $L$ and its terminator~$c$.}
        \label{fig:3dpolar}
        \end{figure}

           Let us have a point light source $L[l_1,\dots,l_n]$ with homogeneous coordinates $\overline{L}(\overline{l_1},\dots,\overline{l_n},\overline{l_0})$ . The terminator $c$ of the hypersurface $\mathcal{S}$ with respect to $L$ is the intersection of the first polar
        \begin{equation}
        \mathcal{S}_L: \overline{\sigma_L}=\overline{L}^T\nabla \overline{\sigma}=0
        \end{equation}
        of the hypersurface $\mathcal{S}$ with the hypersurface $\mathcal{S}$:
        
      \begin{equation}
        \label{eq:contour}
        c: \sigma=0 \wedge \sigma_L=0,          
        \end{equation}
        where $\sigma_L$ is the dehomogenized polynomial $\overline{\sigma_L}$ (Figure~\ref{fig:3dpolar}). 
        

        \subsubsection{Tangent hypercones}
        \label{subsubsection:TangentCones}

        \begin{figure}[!t]
        \centering
           \includegraphics[width=\linewidth]{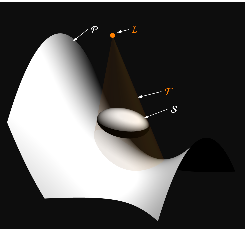}
        \caption{The tangent hypercone $\mathcal{T}$ to a hypersurface $\mathcal{S}$ through a light source $L$ and the shadow cast on $\mathcal{P}$.}
        \label{fig:3dcone}
        \end{figure}
        
        The next step is to find an implicit representation of the tangent hypercones\footnote{For the sake of readability, we use terms cones and hypercones instead of more proper terms conical surfaces, conical hypersurfaces, \dots, over the paper.} from the light source $L$. These are the hypercones through the terminator $c$. Thus, let $Q(q_1,\dots,q_n)$ be a point on $c$. The tangent cone $\mathcal{T}$ is the set of lines $LQ$ for all points $Q$. 
        The line $LQ$ can be (parametrically) represented in a general dimension with the parameter $a\in {\rm I\!R}$ as the set of points $X(x_1,\dots,x_n)$ satisfying the following $n$ equations (one equation for each coordinate):
        \begin{equation}
        aL+(1-a)Q=X.
        \end{equation}
        Thus, the implicit equation of $\mathcal{T}$ is the solution of the system:
        \begin{equation}
        \label{eq:cone}
        \mathcal{T}: \sigma(Q)=0 \wedge \sigma_L(Q)=0 \wedge aL+(1-a)Q-X=0,
        \end{equation}
     where $\sigma(Q)$ and $\sigma_L(Q)$ denote polynomials $\sigma$ and $\sigma_L$ in variables $Q$. 
     By computing the Gr{\" o}bner basis or Dixon resultant of the system 
    \begin{displaymath}       
        \{\sigma(Q),\sigma_L(Q),aL+(1-a)Q-X\}   
    \end{displaymath}
        and eliminating variables $a,q_1,\dots,q_n$, we obtain the polynomial $\theta$ in variables $x_1,\dots,x_n$ representing the tangent hypercone
     \begin{equation}
        \label{eq:conefin}
         \mathcal{T}: \theta=0 \text{~(Figure~\ref{fig:3dcone})}.
     \end{equation}
     
   \subsubsection{Shadow cast on an algebraic hypersurface}
        \label{subsubsection:ShadowCast}

\begin{figure}[!htb]
        \centering
         \includegraphics[width=\linewidth, trim= 5 80 5 40,clip]{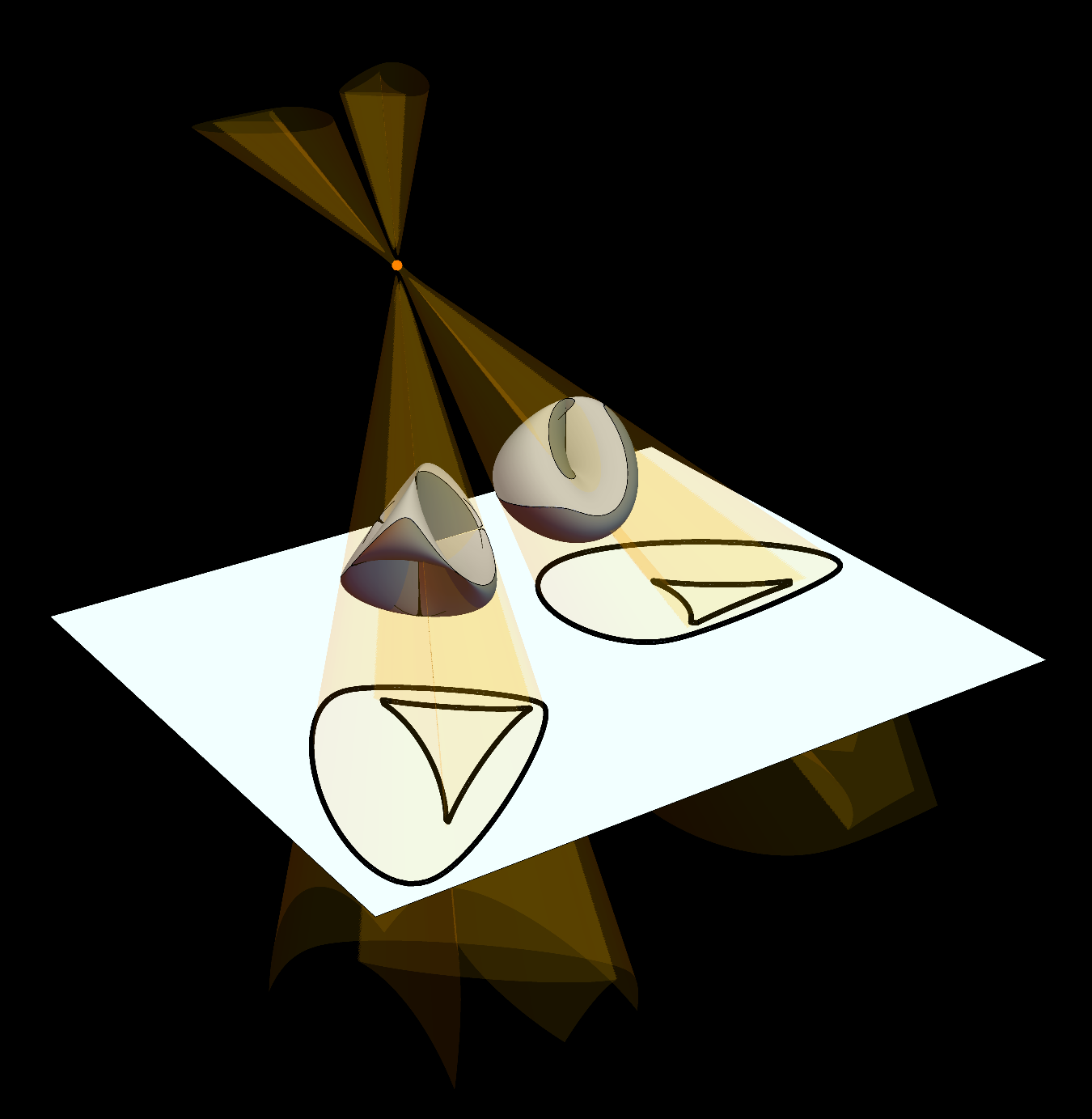}
        \caption{(left) Roman surface: $-2 (x-2) (y-1) z+(x-2)^2 (y-1)^2+(x-2)^2 z^2+(y-1)^2 z^2=0$ and (right) Cross-cap: $\left((x+1)^2+y^2\right) \left((x+1)^2+z^2\right)+(x+1)^2+\frac{y^4}{4}+y^2 z=0$, their parts separated by the first polars, and polar boundaries projected to a plane $z+2=0$. Without omission of the self-shaded parts.}
        \label{fig:sf3d}
    \end{figure}
           Now, we find the shadow cast by $\mathcal{S}$ on itself and on another algebraic hypersurface $\mathcal{P}$ given by a polynomial equation 
        \begin{equation}
        \mathcal{P}:\pi=0.
        \end{equation}

        The boundary of the shadow cast by $\mathcal{S}$ on itself --- terminator, is the intersection of $\mathcal{T}$ with $\mathcal{S}$. 
        The selection of illuminated parts is carried out in the following steps: 
        \begin{enumerate}
            \item  The $n$-space is divided by the first polar to two subspaces ($\sigma_L>0$ or $< 0$), and the illuminated part is in the same subspace as the source of light $L$. The second part is in the shade.\footnote{The method is implemented in 4-D scenes up to this point. The upcoming decomposition to subcones seems, in most cases, computationally unbearable.}
        \end{enumerate}
          Since we cannot always divide the inner and outer parts, the selection fails with non-orientable or self-crossing surfaces (e.g., see Figure~\ref{fig:sf3d}). Hence,  for simplicity, we assume non self-crossing orientable finite hypersurfaces (or at least their terminators are finite).
     \begin{figure}[!htb]
        \centering
           \includegraphics[width=\linewidth]{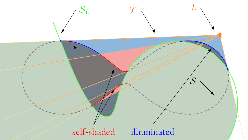}
        \caption{A 2-D situation of a Cassini oval (degree 4) and its self-shading from the point light source. The green curve is the first polar, dividing the plane into two areas. The area that does not contain the light source is excluded. The subcones in 2-D case are plane angles bounded by the rays from the light source. The blue arcs represent illuminated parts, and the red arcs are in the self-shade.}
        \label{fig:cass}
        \end{figure}

    \begin{figure}[!htb]
        \centering
        \includegraphics[width=.49\linewidth,clip]{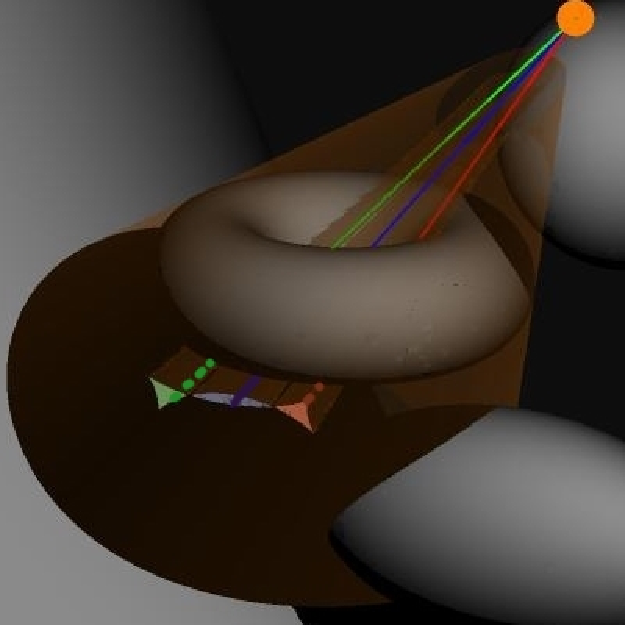}
        \includegraphics[width=.49\linewidth, trim=5 5 5 5,clip]{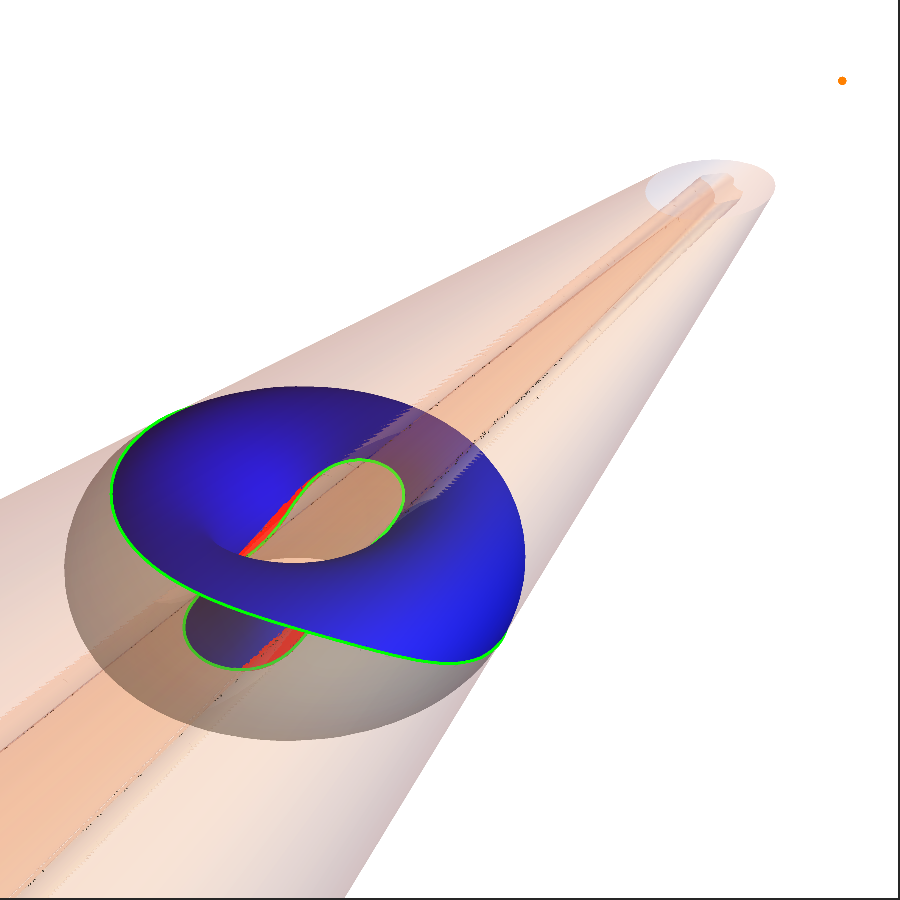}\hfill
        \caption{(left) A shadow of the torus cast on a surface without the self-shaded parts omitted. (right) Decomposition of a torus. The green curve is the terminator. The blue region is illuminated. The red regions are in the accepted subregion divided by the first polar but in the shade of the blue region.}
        \label{fig:tconesproblem}
    \end{figure}
        For polynomials of degrees higher than 2, some regions of the hypersurface $\mathcal{S}$ can still be in their own shade (the tangent hypercone intersects itself), and we have to omit them from the final selection of the illuminated parts (Figures~\ref{fig:cass},\ref{fig:tconesproblem}).
            \begin{enumerate}
             \setcounter{enumi}{1}
            \item Decompose the tangent hypercone into conical subregions (subcones) divided by the terminator.\footnote{We use cylindrical algebraic decomposition. The results of the decomposition are distinct subregions represented by polynomial equations and inequalities (see \cite{Strzebonski2023} for details and implementation).} Eliminate empty subcones.
            \item Select the parts of the hypersurface closest to the light source in non-empty subcones. This is carried out by constructing rays from the light source in each subcone and finding the region with the intersection nearest to the source. The selected nearest parts are illuminated, and the rest is in shade. 
            
        \end{enumerate}

         
         The boundary of the shadow cast by $\mathcal{S}$ on $\mathcal{P}$ is the intersection of its tangent hypercone $\mathcal{T}$ with $\mathcal{P}$. The final shadow contains inner points in the shade, i.e., inside the subcones containing previously selected regions.  

\begin{figure}[htb]
        \centering
     \includegraphics[width=\linewidth]{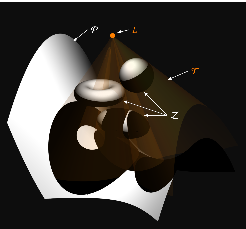}
        \caption{The final visualization of a scene with a sphere, torus, and ellipsoid casting shadows on themselves and on a hyperbolic paraboloid.}
        \label{fig:3Dscene}
        \end{figure}

        If a scene contains more hypersurfaces, some of them might intersect, so we would not be able to distinguish their order with respect to the light source. In such cases, the selection algorithm can be further generalized for a hypersurface $\mathcal{Z}$ as a composition of hypersurfaces $\mathcal{S}_1,\dots,\mathcal{S}_k, k\geq 1$:
        \begin{equation}
        \label{eq:scene}
        \mathcal{Z}: \sigma_1 \dots \sigma_k= 0.
        \end{equation}
        However, this generally works when the light source is outside all composed hypersurfaces. Otherwise, the selection of the illuminated subspaces divided by the first polars must be carried out for some factors separately. For example, in Figure~\ref{fig:3Dscene}, we can define a surface $\mathcal{Z}$ as the composition of a sphere, ellipsoid, and torus. The (infinite) hyperbolic paraboloid would be treated separately.

    \subsection{The principle of a 4-D perspective}
            \begin{figure}[!htb]
        \includegraphics[width=\linewidth]{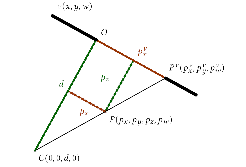}
        \caption{The principle of a 4-D perspective projection from a center $C$ onto a modeling 3-space $\nu(x,y,w)$. 2-D orthogonal projection view on the 4-D perspective shows the homothety between the pre-image in a 4-space and its perspective image in the modeling 3-space from the center of projection.}
        \label{fig:persp2d}
    \end{figure}
        
    Figure \ref{fig:persp2d} shows a 2-D view of the correspondence between the coordinates of a real point $P(p_x,p_y,p_z,p_w)$ and its centrally projected image $P^\nu(p^\nu_x,p^\nu_y,p^\nu_w)$ into the modeling 3-space $\nu(x,y,w)$ (an analogy to a picture plane in a 3-D perspective). Assume an orthogonal coordinate system $(x,y,z,w)$ placed in the point $O$ in $\nu$, with the $z$-axis perpendicular to $\nu$ and the center of projection $C$ in the oriented perspective distance $d$ from $\nu$ such that $CO\perp\nu$. Observing the homothety, the coordinates of the 4-D perspective image $P^\nu$ of $P\neq C$ are 
    \begin{equation}
    \label{eq:persprules}
    (p^\nu_x,p^\nu_y,p^\nu_w)=(d\frac{p_x}{d- p_z},d\frac{p_y}{d-p_z},d\frac{p_w}{d-p_z})
    \end{equation}
    (in the case of $p_z=0$ the coordinates do not change; if $p_z=d$ the image is improper). 

    To find the implicit representation of the 4-D perspective image (3-D occluding contour) of an algebraic surface 
    \begin{equation}
    \mathcal{S}:\sigma=0
    \end{equation}
    (in variables $x,y,z,w$), we use a similar idea as in constructing shadows. Let 
    \begin{equation}
    \mathcal{S}_C:\overline{\sigma_C}=\overline{C}^T\nabla \overline{\sigma}=0,
    \end{equation}
    in homogeneous coordinates, be the first polar of $\mathcal{S}$ with respect to $C$, and with \begin{equation}
        \sigma_C=0
    \end{equation} being its dehomogenized equation. Since we are projecting to a 3-space instead of an arbitrary hypersurface, we can use the derived rules of the 4-D perspective mapping (Equations~\ref{eq:persprules}). In the next step, we set up a system of polynomials prepared for elimination using Gr{\" o}bner basis or Dixon resultant such that the final image will be in the coordinates $(x,y,w)$ of the modeling space. 
    
    Let $Q(q_x,q_y,q_z,q_w)$ be a point on a contour generator 
    \begin{equation}
    \label{eq:cgreal}
    c_C:\sigma=0 \wedge \sigma_C=0.
    \end{equation}
    First, we substitute variables $X(x,y,z,w)$ by $Q(q_x,q_y,q_z,q_w)$ in (\ref{eq:cgreal}) such that $\sigma(X)\rightarrow \sigma(Q)$ and $\sigma_C(X)\rightarrow\sigma_C(Q)$. The rules of the mapping (\ref{eq:persprules}) are represented by the system of linear equations: 

     \begin{equation}
	   \label{eq:polyrules}
	   \begin{split} 
	   x-d \frac{q_x}{d-q_z}=0,\\
	     y-d \frac{q_y}{d-q_z}=0,\\
	     w-d \frac{q_w}{d-q_z}=0.
	   \end{split}
    \end{equation}

   Elimination of $q_x,q_y,q_z$, and $q_w$ from the system of polynomials 
   \begin{displaymath}
        \biggl\{ \sigma(Q),\sigma_C(Q),  x-d \frac{q_x}{d-q_z},  y-d \frac{q_y}{d-q_z},  w-d \frac{q_w}{d-q_z}\biggr\}
        \end{displaymath}
        leads to a polynomial $\sigma^\nu$ in $(x,y,w)$ such that its zero set represents the 4-D perspective image $\mathcal{S}^\nu$ of the surface~$\mathcal{S}$. 

    Assuming a point light source $L$, the terminator $c^\nu$ is, in this case, a 2-surface, obtained by the intersection of $\mathcal{S}$ and $\mathcal{S}_L$, given by Equations~\ref{eq:contour} and projected into modeling space $\nu$ through Equations~\ref{eq:persprules}. The images of terminator 2-surfaces are  derived from the system of polynomials  
    \begin{displaymath}
    \biggl\{ \sigma(Q),\sigma_L(Q),  x-d \frac{q_x}{d-q_z},  y-d \frac{q_y}{d-q_z},  w-d \frac{q_w}{d-q_z}\biggr\},
    \end{displaymath} 
    similarly as above.

    For the sake of better understanding, we also visualize the tangent hypercones $\mathcal{T}$ of $\mathcal{S}$ from $L$, when possible (cf. Subsections~\ref{subsection:4dmoon} and~\ref{subsection:4dring}). The tangent hypercones $\mathcal{T}$ are 3-surfaces given by Equation~\ref{eq:conefin}, and the contours of their images $\mathcal{T}^\nu$ are created by the same procedure as images $\mathcal{S}^\nu$ of $\mathcal{S}$. 
    
    The final shadows are three-dimensional regions bounded by 2-spaces in the 4-space. Conveniently, in a 4-D perspective, we can find implicit equations of the 2-surface boundaries in the modeling 3-space. To do so, we need to find the zero set of the system of polynomials  \begin{displaymath}
    \biggl\{ \sigma(Q),\theta(Q), x-d \frac{q_x}{d-q_z},  y-d \frac{q_y}{d-q_z},  w-d \frac{q_w}{d-q_z}\biggr\},
    \end{displaymath}
    representing the perspective image of the intersection of the surface $\mathcal{S}$ and the tangent cone $\mathcal{T}$.  The last step is to select the regions according to Section~\ref{subsubsection:ShadowCast}. In a 4-D perspective, we only highlight shadow boundaries so that we can see through 3-D images. 
    
\section{Experimental results and technical details}
In this section, we review several examples and comment on technical details. 
\subsection{3-D scene}
\subsubsection{“Implicit Bakery'', Figure~\ref{fig:3Dscene}}
See the text file with equations and video in Appendix~\ref{att:3dbakery}.

The 3-D scenario from Section~\ref{sec:method} shows a surface $\mathcal{Z}$ composed of three implicitly given factor surfaces: sphere $\mathcal{S}_1$, torus $\mathcal{S}_2$, and ellipsoid $\mathcal{S}_3$

\begin{equation}
    	\label{eq:3dsphere} \mathcal{S}_1: (x - 1)^2 + (y + 4)^2 + (z - 5)^2 - 4 = 0, 
\end{equation}
\begin{equation}
\begin{split}
     \label{eq:3dtorus} \mathcal{S}_2: \left((x - 1)^2 + (y - 1)^2 + (z - 2)^2 + 3\right)^2\\ - 16(x - 1)^2 + 16(y - 1)^2=0, 
\end{split}
\end{equation}
\begin{equation}
     \label{eq:3dellipsoid} \mathcal{S}_3: 4(x - 3)^2 + (y + 1)^2 + 4(z + 2)^2 -12=0. 
\end{equation}

We describe constructions of shadows cast between them and their shadows cast on a hyperbolic paraboloid~$\mathcal{P}$
\begin{equation}
\mathcal{P}: 2(y + 3)^2 - 2(x - 5)^2 - 25(z + 7)=0
\label{eq:3dhp}
\end{equation}
from the light source $L\left[-1, -2, 10\right]$.

The polynomial of $\mathcal{Z}$ has degree 8, but we treat its decomposition into factors (factor surfaces $\mathcal{S}_1, \mathcal{S}_2, \mathcal{S}_3)$. In this case, the factor surface of the highest degree, 4, is the torus.

A terminator line of each "factor surface" is its intersection with the first polar with respect to $L$. Explicitly, terminators are the four sets of points that satisfy the following pairs of equations:

\begin{center}
    
  \begin{tabular}{
				@{\extracolsep{\fill}}
				l
                >{$\displaystyle$}c<{$\vphantom{1}$}
				>{$\displaystyle$}r<{$\vphantom{1}$}
				>{\refstepcounter{equation}(\theequation)}r
				@{}
			}
    \toprule
    factor surface &  & first polar & \multicolumn{1}{c}{}\\
    \midrule
    	$\mathcal{S}_1$: Eq. \ref{eq:3dsphere} &  $\wedge$ & $19 + 2 x - 2 y - 5 z=0$ &\label{eq:polar3Dsphere}	\\
     \midrule
     $\mathcal{S}_2$: Eq. \ref{eq:3dtorus} &  $\wedge$ & \begin{tabular}{r}      
            
	          $128 - 30 x + 12 x^2$\\ $+ 2 x^3- 29 y - 10 x y$\\$ + 3 x^2 y + 10 y^2 +  2 x y^2$\\$ + 3 y^3- 104 z+ 8 x z $\\$- 8 x^2 z + 4 y z - 8 y^2 z $\\$+  40 z^2+ 2 x z^2 + 3 y z^2 $\\$ - 8 z^3=0$
		 \end{tabular}& \label{eq:polar3Dtorus}\\
   \midrule
     $\mathcal{S}_3$: Eq. \ref{eq:3dellipsoid} &  $\wedge$ &  $16 x + y - 48 z -131 =0$ & \label{eq:polar3Dellipsoid}\\
     \midrule
     $\mathcal{P}$: Eq. \ref{eq:3dhp}	& $\wedge$ & $ 24 x + 4 y - 25 z-708=0$ & \label{eq:polar3Dhp}\\
    \bottomrule
  \end{tabular}
  \end{center}

If a surface has a degree $n$, its first polar has a degree $(n-1)$, and the degree of the terminator line is due to Bézout's theorem $n(n-1)$. Therefore, the terminator line of the entire scene $\mathcal{Z}$ would have the degree $8\cdot 7 =56$. However, after factorization, the factor surfaces will have degrees $2$ for the quadrics and $12$ for the torus. The same holds for degrees of tangent cones \cite{Bydzovsky1948}. 

Next, we omit the shaded subspaces divided by the first polar. In this case, the points on illuminated regions have non-negative values in the equations of the first polars.

\begin{figure}[!htb]%
	\centering
	 \setlength{\lineskip}{-7pt}
	\subfloat{{\includegraphics[width=0.49\linewidth]{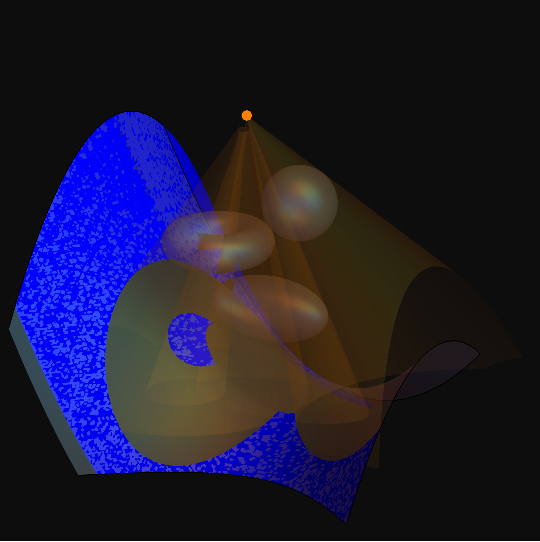} }}%
	\subfloat{{\includegraphics[width=0.49\linewidth]{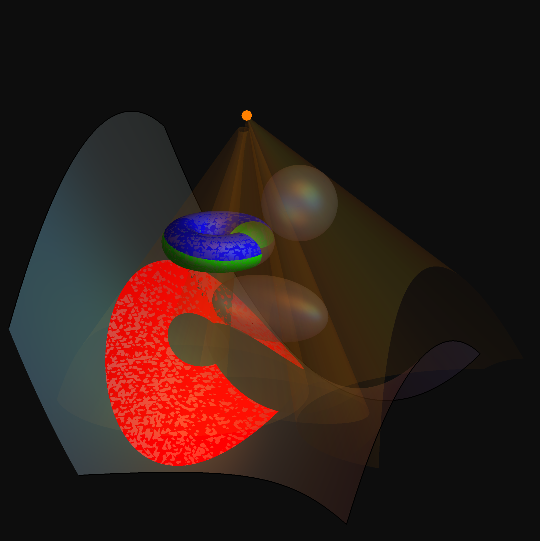} }}%
         \newline
	\subfloat{{\includegraphics[width=0.49\linewidth]{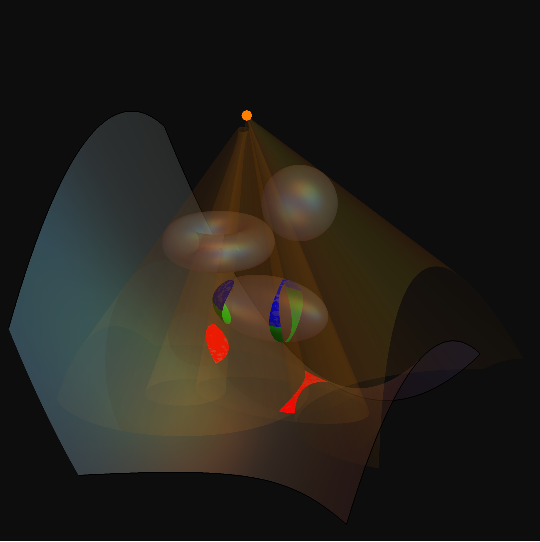} }}%
	\subfloat{{\includegraphics[width=0.49\linewidth]{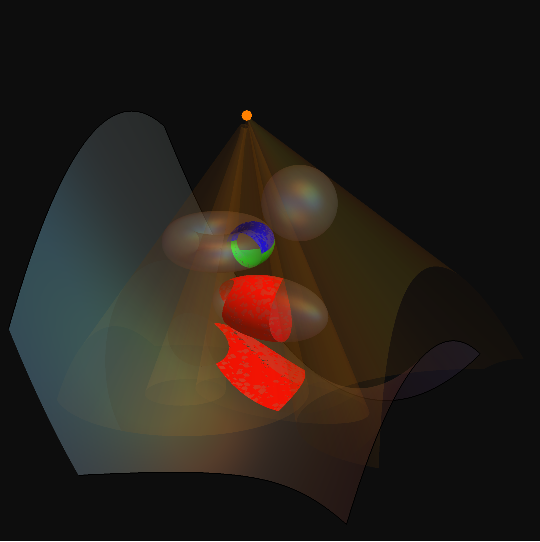} }}%
        \newline	
        \subfloat{{\includegraphics[width=0.49\linewidth]{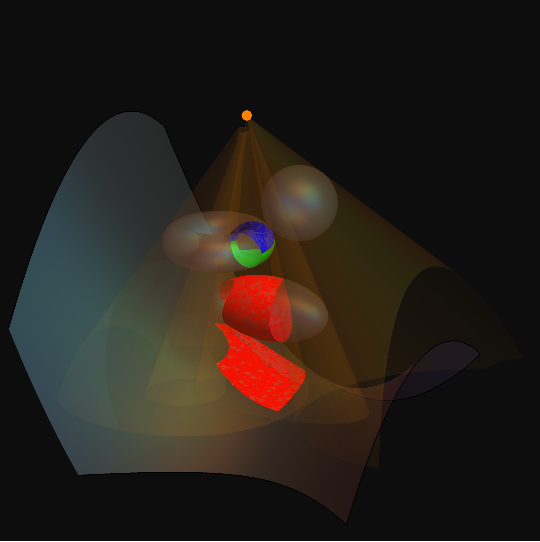} }}%
	\subfloat{{\includegraphics[width=0.49\linewidth]{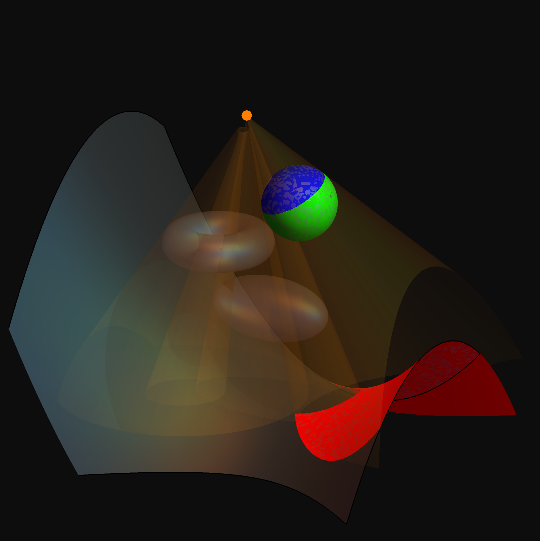} }}%
 \newline
	\caption{A selection of the regions contained in subcones. Each picture shows a different subcone of the scene $\mathcal{Z}$ consisting of a sphere $\mathcal{S}_1$, torus $\mathcal{S}_2$, and ellipsoid $\mathcal{S}_3$. The blue regions are illuminated, the green regions are excluded by the first polars, and the red regions are in the shade of the blue regions with a shorter distance to the light source. The top left figure shows the complement of the union of all subcones in the 3-space; hence it distinguishes the shadow of the scene on the hyperbolic paraboloid $\mathcal{P}$. The rest of the empty regions are not shown.}%
	\label{fig:subcones}%
\end{figure}

Tangent cones are also treated separately for each factor surface. For example, for the torus $\mathcal{S}_2$, we have the following system of polynomial equations:
\begin{equation} 
    \begin{split}
\textrm{Eq.~(\ref{eq:3dtorus})},\\
\textrm{Eq.~(\ref{eq:polar3Dtorus})},\\
aq_1-(a-1)-x=0,\\
aq_2-2(a-1)-y=0,\\
aq_3+10(1-a)-z=0.
\end{split}
\end{equation}
Eliminating $q_1,q_2,q_3$, and $a$ leads to an 8th-degree polynomial of the tangent cone. Similarly, we find the rest of the tangent cones of the surfaces in the scene. 

Next, we decompose the regions bounded by the system of tangent cones into all subcones and obtain five nonempty and three empty intersecting regions in this case (Figure~\ref{fig:subcones}). The subcones are represented by implicit equations and inequalities.

In the last step, we trace subspaces separated by the first polars over each subcone and choose the region nearest to the light source. 

\subsection{4-D scenes}
\subsubsection{Understanding the 4-D visualizations}
Let us give a few remarks on how to understand the 4-D visualizations below: 
\begin{itemize}
    \item The visualizations are 3-D models (occluding contours), and the figures in the paper are only 2-D images of the 3-D scenes.
    \item Standing in a gallery in front of an actual 2-D painting in the 3-D linear perspective, we can find the position of our eye such that the picture makes an illusion of 3-D space. This is unreachable in the 4-D perspective because we cannot leave the 3-D space of the 3-D image (of a 4-D object).
    \item One could think of observing the picture from the inside, i.e., we can reach any location in the 3-D static image. On the contrary, the change of the 4-D eye/camera position would change the contours of the objects.
    \item For better orientation, we always attach an image of a reference hypercube with the center in $[0,0,0,0]$ (Figure~\ref{fig:hypcube}).
      \begin{figure}[!hb]
        \centering
        \includegraphics[width=.25\linewidth]{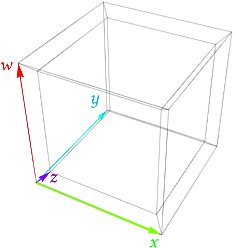}
        \caption{A reference hypercube in the 4-D perspective. The colors of the arrows mark the directions of the coordinate axes: {\color{green}Green $\rightarrow x$}, {\color{cyan}Blue $\rightarrow y$}, {\color{violet} Purple $\rightarrow z$}, {\color{red} Red $\rightarrow w$}.}
        \label{fig:hypcube}
    \end{figure}
  \item The default value of the oriented eye distance in figures is \mbox{$d=-6$}, and its coordinates are $[0,0,d,0]$.
    \item To understand the (3-D) spatial properties of the modeling 3-space, we keep the software lighting properties of the 3-D graphics.  
\end{itemize}

\subsubsection{Technical notes on the implementation of the 4-D perspective}
Let us reflect on some pros and cons of the 4-D perspective:
    \begin{itemize}
        \item[$+$] Intersections of 3-surfaces in a 4-space are \mbox{2-surfaces}; hence, we can visualize them using only one implicit equation in the modeling 3-space. The equation can be obtained directly from the corresponding polynomial system. 
        \item[$-$] 4-D perspective images should include "inner points" of hypersurfaces, but we only show their occluding contours. Otherwise, we would not see images that overlap in the modeling 3-space. An analogous problem occurs in 3-D perspective, where we can imagine a smaller object in front of a bigger object, so the perspective image of the smaller object would lie inside the image of the bigger object. Although we can easily decompose figures in the picture plane, we would not see much in the modeling 3-space. The understanding becomes very unclear with non-closed surfaces, where parts of hypersurfaces might seem to go through their contours (see Section~\ref{subsection:4dmoon}).
        \item[$-$] The algorithm to find and fill the unshaded and shaded parts of surfaces inside tangent cones works well in theory and is successfully applied in 3-D scenes. However, it meets technical difficulties in 4-D. In our experiments, the computation time to find points on the hypersurface (polynomial equation in four variables) and inside the subcone (system of polynomial inequalities in four variables) was beyond reasonable limits. For this purpose, a parametric representation might be more appropriate.  
    \end{itemize}

\subsubsection{4-D scene: “HyperQuadrics'' Figure~\ref{fig:hyperquadricsfull}}
See the text file with equations, code, and video in in Appendix~\ref{att:4dhyperquadrics}.  

\begin{figure}[!htb]
    \centering
    \includegraphics[width=\linewidth]{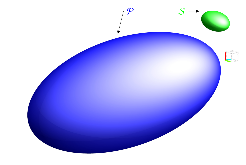}
    \caption{A 4-D scene with a 3-sphere $\mathcal{P}$ and a 3-ellip\-so\-id~$\mathcal{S}$ in a 4-D perspective.}
    \label{fig:hyperquadricspoly}
\end{figure}

In the first 4-D case (Figure~\ref{fig:hyperquadricspoly}), we have a 3-ellipsoid 
\begin{equation}
    \mathcal{S}: \frac{(x + 2)^2}{4} + \frac{(y + 1)^2}{2} + z^2 + (w - 4)^2 - 1=0
\end{equation}
casting a shadow on a 3-sphere 
\begin{equation}
    \mathcal{P}: (x + 5)^2 + (y + 6)^2 + (z - 2)^2 + (w + 3)^2 - 36=0
\end{equation}
from the light source $L=[1, 1, -1.5, 5]$.

In this scene, the given 3-surfaces and their tangent hypercones have degrees 2. The highest degree, 4, has the two-dimensional boundary of the shadow cast by the \mbox{3-ellipsoid} on the 3-sphere. The occluding contours $\mathcal{S}^\nu$ and $\mathcal{P}^\nu$  of $\mathcal{S}$ and $\mathcal{P}$ are ellipsoids, and the same holds for their terminator 2-surfaces $c_\mathcal{S}, c_\mathcal{P}$  and their images $c^\nu_\mathcal{S}, c^\nu_\mathcal{P}$.

\begin{figure}[!b]
    \centering
    \includegraphics[width=.3\linewidth,trim= 0 50 0 0,clip]{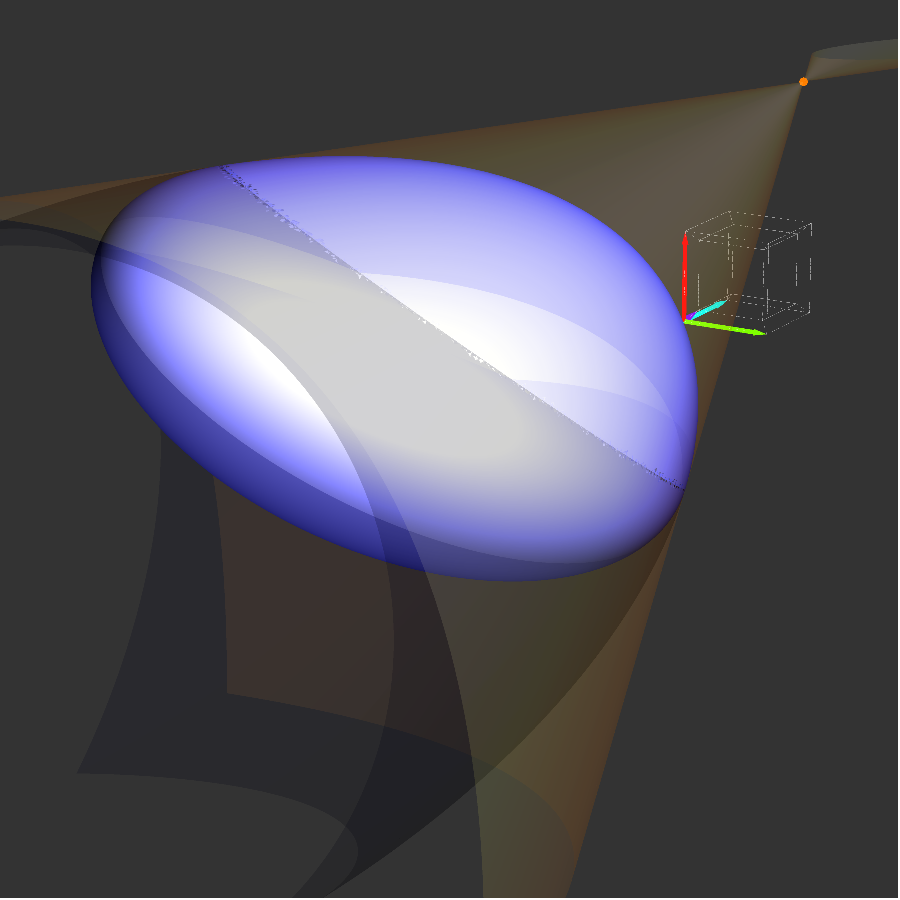}
    \includegraphics[width=.3\linewidth,trim= 0 50 0 0,clip]{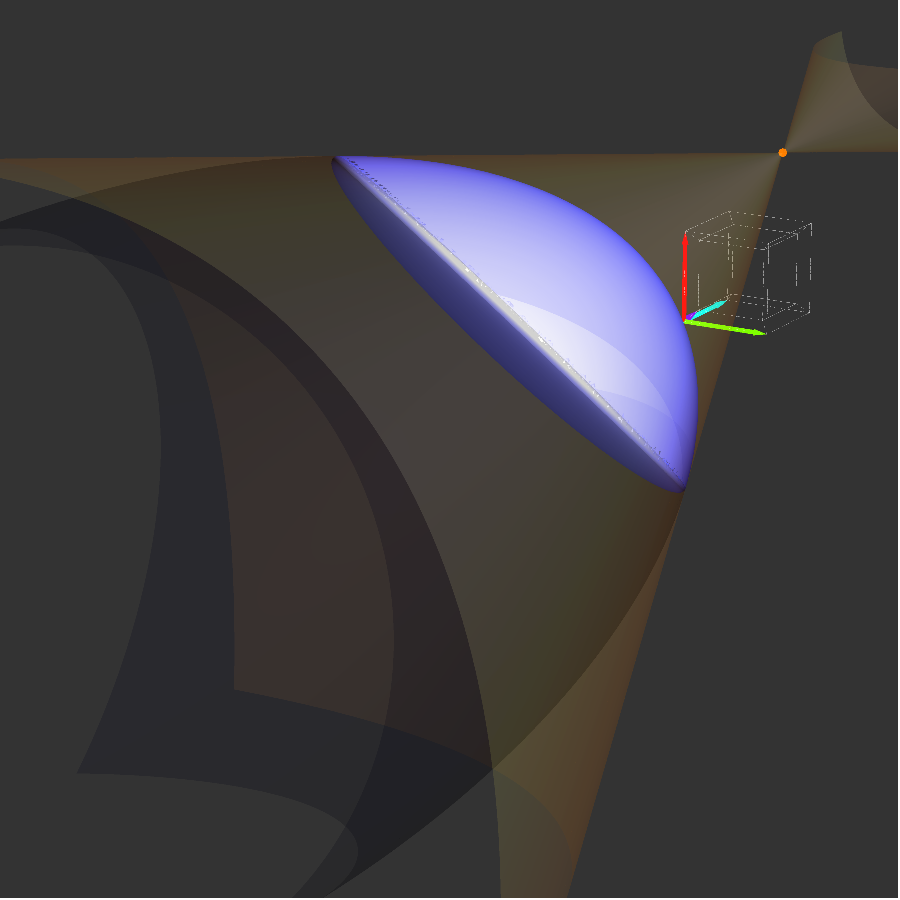}
    \includegraphics[width=.3\linewidth,trim= 0 50 0 0,clip]{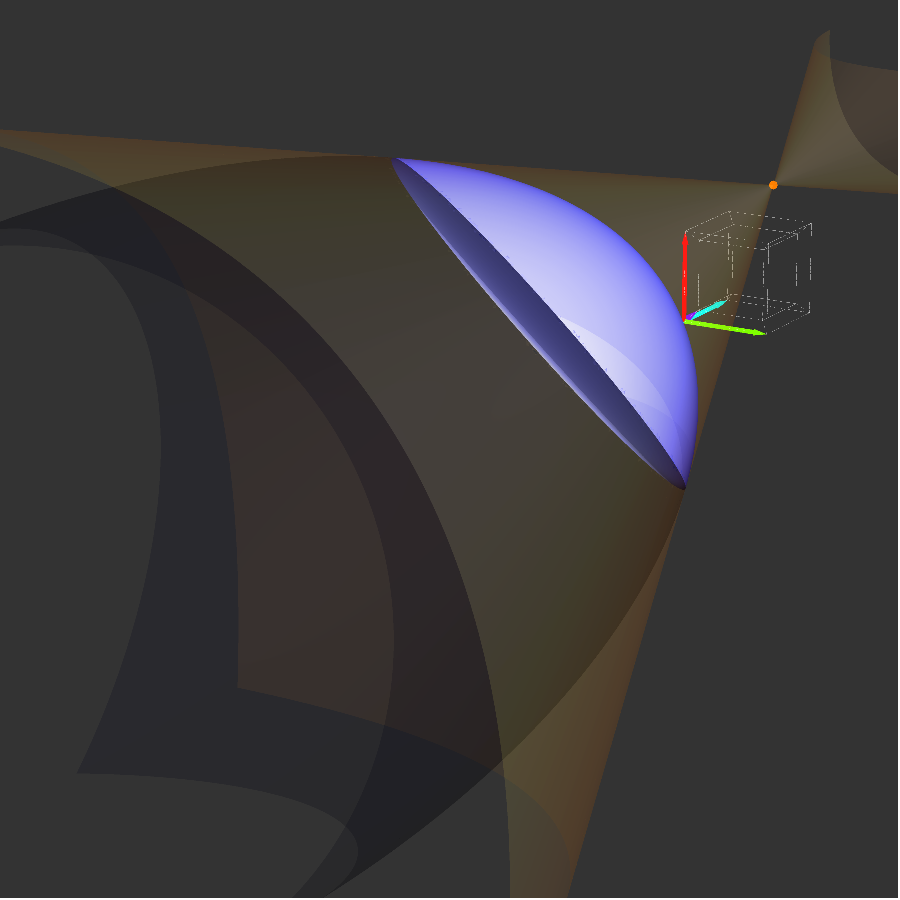}
    \caption{Transitions of the illuminated part of the \mbox{3-sphere}~$\mathcal{P}$ from the point light source $L=[1, 1, z_L, 5]$ moving in the $z$-direction: (left) $z_L=2.5$, (center) $z_L=7.5$, (right) $z_L=12.5$. The point $L$ lies in the first polar 3-space for $z_L=9.5$. The illuminated part is bounded by the occluding contours of the 3-sphere and terminator 2-surface with respect to $L$.}
    \label{fig:4dspheretrans}
\end{figure}

In the case of hyperquadrics, we can easily deduce the transition of the visible illuminated parts with respect to the given 4-D perspective (Figure~\ref{fig:4dspheretrans}). The first polar of the 3-sphere $\mathcal{P}$ with respect to the perspective center $C$ is a 3-space, dividing the 4-space into two half-4-spaces. Thus, we have the following three cases:
\begin{enumerate}
\item The shape of the illuminated part is in a special position when the light source $L$ is in the polar 3-space (with respect to $C$), i.e., the terminator 2-surface of the 3-sphere with respect to $L$ degenerates to an ellipse including inner points. 
\item When $L$ is in the same half-4-space as the center $C$, the visible illuminated part includes the inner points of the terminator 2-surface. 
\item If the light source $L$ and the center $C$ lie in the opposite half-4-spaces, we must exclude the inner points of the terminator 2-surface.
\end{enumerate}

\begin{figure}[!htb]
    \centering
    \includegraphics[width=\linewidth]{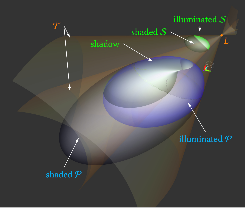}
    \caption{A 4-D scene with the shadow of a 3-ellipsoid~$\mathcal{S}$ on a 3-sphere.}
    \label{fig:hyperquadricsfull}
\end{figure}

\subsubsection{4-D scene: "Full HyperMoon Between HyperMountains", Figure~\ref{fig:moonshadow}}
See the text file with equations and code in Appendix~\ref{att:4dhypermoon}. 

\begin{figure}[!htb]
    \centering
    \includegraphics[width=\linewidth]{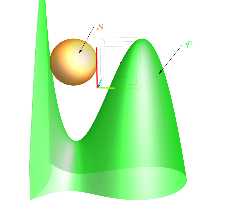}
    \caption{A 4-D scene with a 3-sphere $\mathcal{S}$ and a 3-sur\-face~$\mathcal{P}$ of degree 3 in a 4-D perspective.}
    \label{fig:moon}
\end{figure}

\label{subsection:4dmoon}
The second 4-D scene (Figure~\ref{fig:moon}) shows a more complicated situation.  Let us have a 3-sphere 
\begin{equation}
    \mathcal{S}: w^2+(x+1)^2+y^2+(z+1)^2-\frac{1}{4} =0
\end{equation}
casting a shadow on a 3-surface of degree 3
\begin{equation}
    \mathcal{P}: (x-1) (x+2) x+y^2+z^2+w=0
\end{equation}
from the light source $L=[5, 1,1, 2]$.

The contours of the occlusion are 2-surfaces of degree~2 for $\mathcal{S}^\nu$ and degree 6 for $\mathcal{P}^\nu$. 
In this case, it is hard to perceive the 4-D spatial properties of the scene from the contours. In particular, we cannot intuitively grasp the infinite 3-surface $\mathcal{P}$.   
\begin{figure}[!htb]
    \centering
    \includegraphics[width=\linewidth]{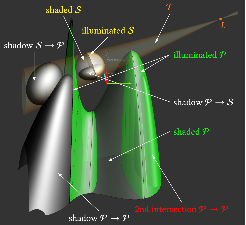}
    \caption{An illumination of a 3-sphere $\mathcal{S}$ and its shadow on a 3-surface~$\mathcal{P}$ of degree 3 from a point light source $L$ in a 4-D perspective. The figure contains the excess intersection of the tangent cone (not visualized) of~$\mathcal{P}$, and the self-shaded region of $\mathcal{P}$ is not excluded from the illuminated part.}
    \label{fig:moonshadow}
\end{figure}

Let us bring more light to this scene. After finding the terminator 2-surfaces and tangent hypercones\footnote{The computations of the projection of the tangent hypercone to $\mathcal{P}$ was terminated after too long ($5-10$h).} to the 3-sur\-faces through the vertex $L$,
we can create shadows between the 3-surfaces. The 3-sphere casts a shadow on $\mathcal{P}$. The contour of the intersection of the tangent hypercone to $\mathcal{S}$ with $\mathcal{P}$ is a 2-surface of degree~$6$. It consists of two disjoint parts, and the part closer to $L$ is omitted in Figure~\ref{fig:moonshadow}. 
The image of the 2-surface boundary of the shadow of $\mathcal{P}$ cast on $\mathcal{S}$ is given by a polynomial of degree 14 
and the contour of the shadow of $\mathcal{P}$ on itself is a surface of degree 18.

\subsubsection{4-D scene: “HyperRing", Figure~\ref{fig:ringfull}}
See the text file with equations, code, and video in Appendix~\ref{att:4dhyperring}. 

\label{subsection:4dring}
\begin{figure}[!htb]
   \centering
  \includegraphics[width=\linewidth, trim= 0 100 0 0, clip]{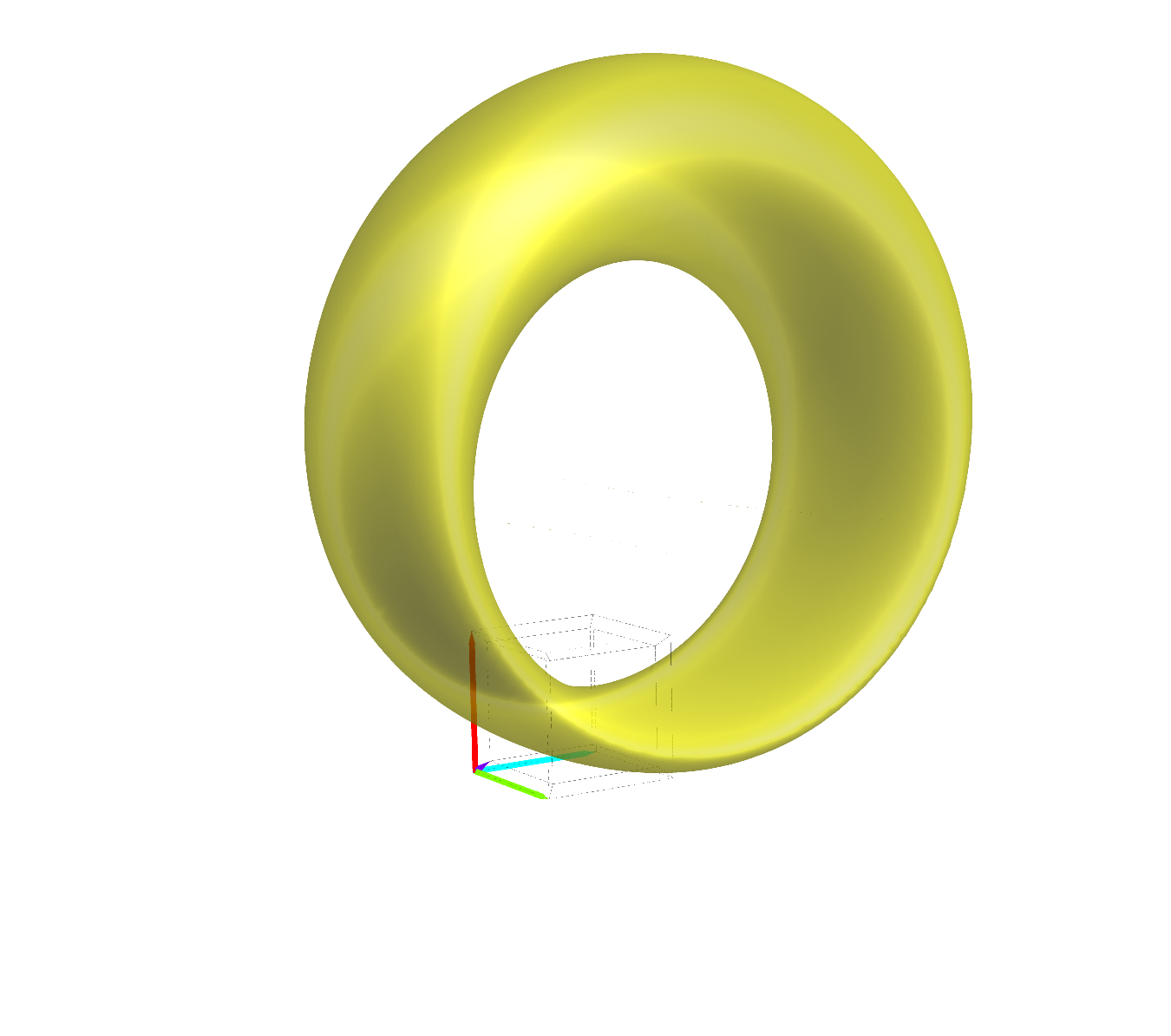}
    \caption{A surface $\mathcal{S}$ of degree 4 in a 4-D perspective.}
\label{fig:ring}
\end{figure}
The last 3-surface (Figure~\ref{fig:moon}) is given by a polynomial of degree~$4$:
\begin{equation}
    \mathcal{S}: (x-1)^2 + ((w-2)^2 + y^2-4)^2 + z^2-1=0.
\end{equation}
The first polar $\mathcal{S}_L$ with respect to the light source $L[0,2,-2,4]$ is:
\begin{equation}
    \begin{split}
    \mathcal{S}_L: 4 w^3 - 3 x + x^2 + 4 w^2 (-8 + y) - 16 y^2\\ + 4 y^3 +  4 w (16 - 4 y + y^2) - 2 z + z^2 =0.
    \end{split}
\end{equation}
Let us have a 3-space 
\begin{equation}
    \mathcal{P}: w+2=0.
\end{equation}  
The occluding contour $\mathcal{S}^\nu$ is after elimination given by a polynomial of degree $8$ in 72 terms (in variables $x,y,z$). 
After elimination, the terminator surface $c^\nu$ is given by a polynomial of degree $8$ in 146 terms (in variables $x,y,z$). 
The tangent cone $\mathcal{T}$ generated by the terminator $c$ (intersection of $\mathcal{S}$ and $\mathcal{S}_L$) is after elimination given by a polynomial of degree $8$ in 483 terms (in variables $x,y,z,w$). 
The occluding contour of the shadow of $\mathcal{S}$ on $\mathcal{P}$ is a 2-surface of degree 8 (Figure~\ref{fig:ringfull}).

\begin{figure}
    \centering
 \includegraphics[width=\linewidth]{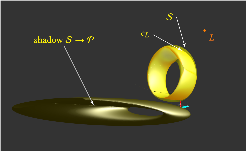}
    \caption{The shadow cast by a 3-surface $\mathcal{S}$ on a \mbox{3-space}~$\mathcal{P}$. Since the terminator overlaps itself, we cannot properly distinguish the illuminated visible part of~$\mathcal{S}$.}
    \label{fig:ringfull}
\end{figure}

\section{Discussion and future work}
Throughout the paper, we tried to bring the "most universal" solution to visualizing shadows of algebraic hypersurfaces. In this sense, the method presented in Section~\ref{subsection:shadow} works in a general dimension; the key point is to construct tangent cones with Equations~\ref{eq:cone}. Illumination of a scene, both in 3-D and 4-D, is a complex process, and we should always consider the properties and positions of objects in the scene. The critical features of objects in our approach are the finiteness, orientability, and degree of a hypersurface and its terminator. Additionally, the position of the point light source and the projection center is crucial to define the inside and outside of a hypersurface or its part and, consequently, its visibility. The degree of a hypersurface opens up problems of computational complexity. Evidently, the appropriate setting and the choice of elimination method (Gr{\" o}bner basis or the Dixon resultant) plays a crucial role in the computation time for higher-degree polynomial systems.
The computational complexity also depends on the form of the polynomial, i.e., transformations of the hypersurface. The computations with hypersurfaces of degree four and higher easily failed after minor adjustments in our experiments. Furthermore, the diversion between dimensions became important in the final visualizations of the inner points of shadows. While in 3-D, we could cover the regions with a satisfying number of points, our solution in 4-D was defeated by time to solve a system of equations and inequalities in four variables. We dropped this case because filling the 3-dimensional volumes would not clarify our visualizations in any way. 

The question of perceiving four- and more-dimensional spaces is very challenging. We only added one more "perspective" open for further investigation. One of the possibilities is to study the properties of 3-surfaces through their projections. Since we developed our method on implicit surfaces, it is convenient for mathematical visualization. Second, we can pursue more "natural" illumination details, including 4-dimensional light intensity, specularity, reflections, etc. Polynomial systems that create tangent cones and occluding contours always contain first-degree polynomials. Hence, there might be possibilities for improving the algorithms tailored to our situation and consequently reducing the computation time. Last but not least, illuminated four-dimensional scenes might be used to study, understand, and train four-dimensional spatial ability (if possible).  

\section{Conclusion}
This paper focuses on a four-dimensional visualization based on implicit representations of hypersurfaces. We have described a general method to find shadow boundaries in an arbitrary dimension and applied it in a three- and four-dimensional space. Furthermore, we have designed a system of polynomial equations to construct occluding contours of hypersurfaces in a 4-D perspective. The method was presented on a composed 3-D scene and three 4-D cases with gradual complexity. Our experimental journey was extensively commented on throughout the article.

\section{Acknowledgments}
Jakub \v{R}ada was supported by the grant SVV 260580. We thank the community at mathematica.stackexchange.com for support with minor technical details in {\sl Mathematica}. We also thank Daniel Lichtblau from {\sl Wolfram Research} for his remarks leading to a significant acceleration of our computations in {\sl Mathematica} and an update of the Dixon resultant implementation throughout our research. 

\printbibliography
\clearpage
\begin{appendices}

The attachments are available in the online repository \url{https://github.com/mbzamboj/4-D-shadows/v2} .

\section{Attachment 1} 
\label{att:3dbakery}
\begin{itemize}
    \item Equations: attachment01-3D-bakery.rtf
    \item Video: \url{https://youtu.be/tKqZn7tzAQE}
\end{itemize}
\section{Attachment 2} 
\label{att:4dhyperquadrics}
\begin{itemize}
    \item Equations: attachment02-4D-hyperquadrics.rtf
    \item Code: code02-4D-hyperquadrics.nb
    \item Video: \url{https://youtu.be/01kYwSblPEY}
\end{itemize}
\section{Attachment 3} 
\label{att:4dhypermoon}
\begin{itemize}
    \item Equations: attachment03-4D-hypermoon.rtf
    \item Code: code03-4D-hypermoon.nb

\end{itemize}
\section{Attachment 4} 
\label{att:4dhyperring}
\begin{itemize}
    \item Equations: attachment04-4D-hyperring.rtf
    \item Code: code04-4D-hyperring.nb
    \item Video: \url{https://youtu.be/6ZdwJ-P18Gw}
\end{itemize}
\section{Attachment 5}
Table~\ref{tab1} (Computation times: 4-D Scene: "Full HyperMoon Between HyperMountains" \ref{subsection:4dmoon})\\\\
Table~\ref{tab2} (Computation times: 4-D Scene: "HyperRing" \ref{subsection:4dring})\\\\
Notation:\\
deg -- degrees of polynomials\\
WM-GB -- Wolfram Mathematica implementation of the function GroebnerBasis with the attributes:
\begin{itemize}
\item[] LO: default lexicographic monomial order
\item[] MO (BEO):\\MonomialOrder $\rightarrow$ EliminationOrder,\\Method $\rightarrow$ "Buchberger" 
\item[] MO (EE):\\MonomialOrder $\rightarrow$ EliminationOrder,\\Method $\rightarrow$ \{"GroebnerWalk",\\ "EarlyElimination$\rightarrow$True"\}
\item[] MO (GWEE):\\MonomialOrder $\rightarrow$ EliminationOrder,\\Method $\rightarrow$ \{"GroebnerWalk",\\ "EarlyElimination$\rightarrow$True"\}
\end{itemize}
WM-Dix -- function DixonResultant in Wolfram Mathematica\\
Fer-Dix-EDF -- Dixon resultant implementation in Fermat\\
T -- terminated after 5 hours or more\\
F -- failed

\begin{table*}[!t]
\centering
\caption{Computation times: 4-D Scene: "Full HyperMoon Between HyperMountains" \ref{subsection:4dmoon}}
\label{tab1}
\begin{tabular}{| l | c || c | c | c | c |}
\hline
object & deg & WM-GB LO & WM-GB MO & WM-Dix & Fer-Dix-EDF\\
\hline\hline 
occ. cont. $\mathcal{S}^\nu,\mathcal{P}^\nu$ & 3, 6 & 0.03s &0.01s (EO) & 0.14s & ---\\
\hline 
terminators $c^\nu_\mathcal{S}, c^\nu_\mathcal{P}$ &  2, 6 & 0.26s & 0.01s (EO) & 0.19s & ---\\
\hline 
tang. hypcon. $\tau_S, \tau_P$ & 2, 6 & 0.1s & 0.03s (GWEE) & F & ---\\
\hline 
shadow $\mathcal{S}\rightarrow\mathcal{P}$ & 6 & 25s & 0.03s (GWEE) & 0.33s & ---\\
\hline 
shadow $\mathcal{P}\rightarrow\mathcal{P}$ & 18 & T & F & 11s & 26s\\
\hline 
shadow $\mathcal{P}\rightarrow\mathcal{S}$ & 14 & T &0.24s (BEO) & 15966s & 79s\\
\hline 
\end{tabular}

\end{table*}

\begin{table*}[!htb]
\centering
\caption{Computation times: 4-D Scene: "HyperRing" \ref{subsection:4dring}}
\label{tab2}
\begin{tabular}{| l | c || c | c | c | c |}
\hline
object & deg & WM-GB LO & WM-GB MO & WM-Dix & Fer-Dix-EDF\\
\hline\hline 
occ. cont. $\mathcal{S}^\nu$ & 8 & 1.2s & 0.02s (EO) & 0.27s & ---\\
\hline 
terminator $c^\nu_\mathcal{S}$ & 8 & 516s & 0.06s (EO) & 0.40s & ---\\
\hline 
tang. hypcon. $\mathcal{T}$ & 8 & 2701s & 0.20s (EO) & 75s & 3.08s\\
\hline 
shadow $\mathcal{S}\rightarrow\mathcal{P}$ & 8 & 0.06s & 0.05s (BEO) & 8.6s & ---\\
\hline 
\end{tabular}
\end{table*}

\end{appendices}

\end{document}